# Consistency of a needlet spectral estimator on the sphere

**Gilles Faÿ**

*Laboratoire Paul Painlevé, Université Lille 1,*
*and Laboratoire AstroParticule et Cosmologie (APC), Université Paris Diderot - Paris 7,*
*e-mail:* `gilles.fay@univ-lille1.fr`

**Frédéric Guilloux**[*]

*MODAL'X, Université Paris Ouest – Nanterre-la-Défense,*
*Laboratoire de Probabilités et Modèles Aléatoires (LPMA), Université Paris Diderot - Paris 7*
*and Laboratoire AstroParticule et Cosmologie (APC), Université Paris Diderot - Paris 7,*
*e-mail:* `guilloux@math.jussieu.fr`

**Abstract:** The angular power spectrum of a stationary random field on the sphere is estimated from the needlet coefficients of a single realization, observed with increasingly fine resolution. The estimator we consider is similar to the one recently used in practice by (Faÿ et al. 2008) to estimate the power spectrum of the Cosmic Microwave Background. The consistency of the estimator, in the asymptotics of high frequencies, is proved for a model with a stationary Gaussian field corrupted by heteroscedastic noise and missing data.

## Contents

[*]corresponding author.

## 1. Introduction

In many application domains (geophysics, cosmology, hydrodynamics, computer vision, etc.), data are defined on the sphere. If the data fit the model of a stationary stochastic field, their second order characteristics, summarized by the angular power spectrum, is of great importance. It contains all the distribution information in the case of a Gaussian stationary process. Ordinary spherical harmonic transform (SHT), the equivalent of the Fourier Series on the circle, provides a fast and efficient method for spectral estimation in the idealistic case of a fully and perfectly observed sphere.

However, rarely the data are available on the whole sphere. Often it is observed under non-stationary contaminants. This is the case for the cosmic microwave background (CMB) which is a major motivation for this work. For those reasons, during the past decade, *localized* analysis for spherical data has motivated many developments; see [14, 10, 20, 6] and the references therein.

The wavelets provide a powerful framework for dealing with non-stationarities. A recent construction of wavelet frames by [16, 17] has proved efficient to analyze stationary spherical processes, thanks to their good localization property. In the time series literature, wavelets are used for spectral estimation whether in a semi-parametric (see *e.g.* [19]) or a non-parametric [3] context. Our observation model, in addition to being spherical, has the particularity of presenting quite general non stationarity (in the structure of the noise) and we failed to find any reference on the subject even for processes living on Euclidean spaces.

In this paper, we establish the consistency of a spectral estimator constructed on the needlets coefficients in high-frequency asymptotics.

The paper is organized as follows. In Section 2 we present the model, including assumptions on the way the process is sampled. In Section 3 we define the needlet spectral estimators and state the consistency results that hold true under realistic conditions. The finite sample behavior of the estimator is explored by numerical simulations (Section 4). Our conclusion is given in Section 5 and Proofs are postponed to Section 6.

## 2. Model and settings

### *2.1. Gaussian stationary spherical fields*

Consider the unit sphere $\mathbb{S}$ in $\mathbb{R}^3$ with generic element $\xi$. The geodesic distance is given by $d(\xi, \xi') \stackrel{\text{def}}{=} \arccos(\xi \cdot \xi')$ where $\xi \cdot \xi'$ denotes the usual dot product between $\xi$ and $\xi'$ (considered as vectors in $\mathbb{R}^3$). The uniform measure $d\xi$ is the unique positive measure on $\mathbb{S}$ which is invariant by rotation, with total mass $4\pi$. Let $\mathbb{H} \stackrel{\text{def}}{=} \mathbb{L}^2(\mathbb{S}, d\xi)$ be the Hilbert space of complex-valued square integrable functions. We have the following decomposition: $\mathbb{H} = \bigoplus_{\ell=0}^{\infty} \mathbb{H}_\ell$ where $\mathbb{H}_\ell$ is the vector space of spherical harmonics of degree $\ell$, *i.e.* restrictions to the sphere of homogeneous polynomials of degree $\ell$ in $\mathbb{R}^3$ which are harmonic (or, equivalently, the restriction of which are eigenvectors of the spherical Laplacian with eigenvalues $\ell(\ell+1)$). The usual spherical harmonics $\mathcal{Y}_{\ell,m}(\xi)$, $-\ell \leq m \leq \ell$, constitute an orthonormal basis of $\mathbb{H}_\ell$. Therefore, the set of all spherical harmonics $\mathcal{Y}_{\ell,m}$, $\ell \geq 0$, $-\ell \leq m \leq \ell$, is an orthonormal basis of $\mathbb{H}$.

In this paper, we shall be concerned with a zero-mean, mean square continuous and real-valued random field $X(\xi)$. We shall assume that $X$ is second-order stationary, that is $\mathbf{E}\left[X(\rho\xi)X(\rho\xi')\right] = \mathbf{E}\left[X(\xi)X(\xi')\right]$ for all $\rho \in SO(3)$. Then the spherical harmonics coefficients of $X$, $a_{\ell,m} \stackrel{\text{def}}{=} \langle X, \mathcal{Y}_{\ell,m} \rangle_{\mathbb{H}}$, are square integrable random variables which verify $\mathbf{E}[a_{\ell,m} a^*_{\ell',m'}] = \delta_{\ell,\ell'} \delta_{m,m'} C_\ell$ for $m, m' \geq 0$ and $a_{l,-m} = a^*_{\ell,m}$. The inverse spherical harmonics transform reads: $X(\xi) = \sum_{\ell \geq 0} \sum_{m=-\ell}^{\ell} a_{\ell,m} \mathcal{Y}_{\ell,m}(\xi)$. The last equality holds in the sense that $\mathbf{E}\left|\int_{\mathbb{S}} X(\xi) - \sum_{\ell=0}^{L} \sum_{m=-\ell}^{\ell} a_{\ell,m} \mathcal{Y}_{\ell,m}(\xi) d\xi\right|^2 \xrightarrow[L\to\infty]{} 0$. The sequence $(C_\ell)_{\ell \geq 0}$ is called the (angular) power spectrum of $X$. Let $L_\ell$ denote the Legendre polynomial of degree $\ell$ normalized by $L_\ell(1) = \frac{2\ell+1}{4\pi}$. The angular power spectrum is linked to the angular correlation of $X$ by the relation $\sum_{\ell \geq 0} C_\ell L_\ell(\cos\theta) = \mathbf{E}[X(\xi)X(\xi')]$ for all pairs of points such that $d(\xi, \xi') = \theta$. The square integrability of $X$ is equivalent to the condition $\sum_{\ell \geq 0} (2\ell+1)C_\ell < \infty$.

We shall also assume that $X$ is Gaussian. This additional assumption is known to be true if and only if the coefficients $a_{\ell,m}$, $\ell \geq 0$, $m \geq 0$ are independent (see [2] for the "only if" part). As mentioned in the Introduction, the finite-dimensional distributions of a Gaussian stationary field are entirely determined by the second-order characteristics, that is by the angular power spectrum of the field. For instance, the second-order stationarity is equivalent, under Gaussian assumption, to the strict stationarity, *i.e.* for all $\rho \in SO(3)$ and $\xi_1, \ldots, \xi_n \in \mathbb{S}$ the two vectors $(X(\rho\xi_1), \ldots, X(\rho\xi_n))$ and $(X(\xi_1), \ldots, X(\xi_n))$ have the same distribution.

### 2.2. Sampling on the sphere

In any real-life situation, only discretized versions of $X$ are available, and consequently spherical harmonic coefficients are exactly computable only if $X$ is $L$-band-limited, that is all the $a_{\ell,m} = 0$, $\ell > L$, for some $L$ which depends essentially on the number of observed points. The discretization of the sphere and achievement of cubature formulas for geodetic functions is a non-trivial task. During the last decade it was shown ([15, 17]) that there exists a constant $c_0 > 0$ such that for all $L \in \mathbb{N}^*$ there exists a set $(\xi_k, \lambda_k)_{k \in \{1,\ldots,N\}} \in (\mathbb{S} \times \mathbb{R}_+^*)^N$ of cubature points and weights (referred to as a *pixelization of order $L$*) with the following properties.

$$\text{For all } f \in \bigoplus_{\ell=0}^{L} \mathbb{H}_\ell, \quad \int_{\mathbb{S}} f(\xi) d\xi = \sum_{k=1}^{N} \lambda_k f(\xi_k) \quad \text{(cubature formula).} \tag{1a}$$

$$c_0^{-1} L^2 \leq N \leq c_0 L^2. \tag{1b}$$

$$c_0^{-1} L^{-2} \leq \min_{1 \leq k \leq N} \lambda_k \leq \max_{1 \leq k \leq N} \lambda_k \leq c_0 L^{-2}. \tag{1c}$$

$$c_0^{-1} L^{-1} \leq \sup_{\xi \in \mathbb{S}} d(\xi, \{\xi_k\}_{k \in \{1,\ldots,N\}}) \leq c_0 L^{-1}. \tag{1d}$$

$$c_0^{-1} L^{-1} \leq \min_{1 \leq k < k' \leq N} d(\xi_k, \xi_{k'}) \leq c_0 L^{-1}. \tag{1e}$$

The following two lemmas derive straightforwardly from these pixelization properties. The first one is proved in [1] and the second one follows from a simple covering argument.

**Lemma 1.** *There exists a constant $c > 0$ such that for all $L \in \mathbb{N}^*$, $\xi \in \mathbb{S}$ and for all pixelization of order $L$ we have*

$$\sum_{k=1}^{N} \frac{1}{(1 + Ld(\xi, \xi_k))^M} \leq c \quad.$$

**Lemma 2.** *There exists a constant $c > 0$ such that for all $L \in \mathbb{N}^*$, $\xi \in \mathbb{S}$ and $\delta > c_0 L^{-1}$ and for all pixelization of order $L$ we have*

$$c^{-1} \delta^2 L^2 \leq \text{Card}\left\{k \in \{1, \ldots, N\} : d(\xi, \xi_k) \leq \delta\right\} \leq c \delta^2 L^2 \quad.$$

### 2.3. Observation model

We are now in position to give a description of our statistical model. Assume that we observe a noisy and sampled version of low-passed $X$ at successive *scales* $j$, with some missing (or attenuated) data. More precisely, for every $j \in \mathbb{N}$, given some $L_j \in \mathbb{N}^*$,

- let $(\xi_{j,k}, \lambda_{j,k})_{k \in \{1,\ldots,N_j\}}$ be a pixelization of order $4L_j$ ;
- let $W_j : k \in \{1, \ldots, N_j\} \mapsto W_{j,k} \in [0,1]$ and $\sigma_j : k \in \{1, \ldots, N_j\} \mapsto \sigma_{j,k} \in \mathbb{R}_+$ be deterministic, known, functions ;
- and let $\mathbf{B}_j : \ell \in \mathbb{N} \mapsto \mathbf{B}_{j,\ell} \in \mathbb{R}$ such that $\mathbf{B}_{j,\ell} = 1$ if $\ell \leq L_j$ and $\mathbf{B}_{j,\ell} = 0$ if $\ell > 2L_j$.

We observe

$$Y_{j,k} \stackrel{\text{def}}{=} W_{j,k} \left[X_j(\xi_{j,k}) + Z_{j,k}\right], \quad j \in \mathbb{N}, \ k \in \{1, \ldots, N_j\} \tag{2}$$

where $Z_{j,k} \stackrel{\text{def}}{=} \sigma_{j,k} U_{j,k}$ and $U_{j,k}$, $j \in \mathbb{N}$, $k \in \{1, \ldots, N_j\}$ is a triangular array of standard and independent Gaussian random variables and independent of the process $X$. The process $X_j$ is defined by

$$X_j \stackrel{\text{def}}{=} \sum_{\ell \geq 0} \sum_{m=-\ell}^{\ell} \mathbf{B}_{j,\ell} a_{\ell,m} \mathcal{Y}_{\ell,m} \quad , \quad a_{\ell,m} = \langle X, \mathcal{Y}_{\ell,m} \rangle_{\mathbb{H}} \quad . \tag{3}$$

In any CMB experiment, some smoothing is induced by the instrumental beam. Eq (3) is a idealistic version of this low-pass operation.

Without loss of generality, $j \mapsto L_j$ is supposed non-decreasing. In the following, we call $W_j$ the *mask*. The particular case of $W_j$ taking its values in $\{0;1\}$ corresponds to missing data.

In other words, a single realization of $X$ is considered, but independent and noisy measures with an increasing spatial resolution are available. This corresponds, for instance, to the observation model of the CMB. The latter is modeled by astrophysicists as the single realization of a stationary Gaussian process. Its observation is achieved by more and more precise instruments, involving their own observation noise, sky coverage and instrumental beam. Full sky map of moderate resolutions (*e.g.* maps provided by the WMAP collaboration [4]) and observations of small and clean patches of the sky at very high resolution (*e.g.* maps from ACBAR experiment [18]) are available simultaneously. Cosmologists aggregate information for those maps to give a large band estimation of the angular power spectrum.

### 2.4. Needlets and statistical properties of needlet coefficients

#### 2.4.1. General framework

The needlets are second-generation wavelet frames which were introduced by [17]. Let us recall below their definition and first properties, the proofs of which are either referred to existing literature or postponed to Section 6.

Start from the fact that the orthogonal projection on $\mathbb{H}_\ell$ has a kernel involving Legendre polynomials, namely

$$\forall f \in \mathbb{H}, \ (\Pi_{\mathbb{H}_\ell} f)(\xi) \stackrel{\text{def}}{=} \sum_{m=-\ell}^{\ell} \langle f, \mathcal{Y}_{\ell,m} \rangle_{\mathbb{H}} \mathcal{Y}_{\ell,m}(\xi) = \int_{\mathbb{S}} L_\ell(\xi \cdot \xi') f(\xi') d\xi'.$$

Instead of considering single frequencies $\ell$, we shall combine them within frequency bands. For this purpose, define a sequence of functions $b_j : \ell \in \mathbb{N} \mapsto b_{j,\ell} \in \mathbb{R}$, $j \in \mathbb{N}$, called (frequency) *window functions*. The window $b_j$ is supposed to be supported in $[0, L_j^{(b)}]$ for some $L_j^{(b)} \in \mathbb{N}$. The kernels $\Psi_j : (\xi, \xi') \mapsto \sum_{\ell \geq 0} b_{j,\ell} L_\ell(\xi \cdot \xi')$ and $\Lambda_j : (\xi, \xi') \mapsto \sum_{\ell \geq 0} (b_{j,\ell})^2 L_\ell(\xi \cdot \xi')$ have the two following obvious properties. First, for all $f \in \mathbb{H}$, $f(\xi) = \sum_{j \in \mathbb{N}} \int_{\mathbb{S}} \Lambda_j(\xi, \xi') f(\xi') d\xi'$. Second, $\Lambda_j(\xi, \xi') = \int_{\mathbb{S}} \Psi_j(\xi, \xi'') \Psi_j(\xi'', \xi') d\xi''$.

The discretization of the above kernels leads to the following spherical functions called *needlets*. For each scale $j \in \mathbb{N}$, given a pixelization $(\xi_{j,k}, \lambda_{j,k})_{k \in \{1,\ldots,N_j\}}$ of order at least $2L_j^{(b)}$, define

$$\psi_{j,k}(\xi) \stackrel{\text{def}}{=} \sqrt{\lambda_{j,k}} \sum_{\ell \geq 0} b_{j,\ell} L_\ell(\xi \cdot \xi_k) \quad .$$

The needlets $\psi_{j,k}$, $j \in \mathbb{N}$, $k \in \{1, \ldots, N_j\}$ constitute a tight frame of $\mathbb{H}$ [17, 9] if for all $\ell \in \mathbb{N}$, $\sum_{j \in \mathbb{N}} (b_{j,\ell})^2 = 1$. For any (possibly random) function $f$ in $\mathbb{H}$, the coefficients $\langle f, \psi_{j,k} \rangle_{\mathbb{H}}$ are renormalized for sake of notational simplicity: we shall handle the *needlet coefficients*

$$\gamma_{j,k} \stackrel{\text{def}}{=} (\lambda_{j,k})^{-1/2} \langle f, \psi_{j,k} \rangle_{\mathbb{H}} \quad .$$

If $f$ is $L_J^{(b)}$-band-limited, one can compute practically the coefficients $\gamma_{j,k}$, $j \leq J$, in the spherical harmonics domain, from the values of $f$ on the cubature points, as made explicit by the following diagram.
$$(f(\xi_k))_{1\leq k\leq N} \xrightarrow{\text{SHT}} (a_{\ell,m})_{\ell\leq L_J^{(b)}} \xrightarrow{\times} (b_{j,\ell}a_{\ell,m})_{\ell\leq L_j} \xrightarrow{\text{SHT}^{-1}} (\gamma_{j,k})_{1\leq k\leq N_j} \qquad (4)$$

The initial pixelization $(\xi_k, \lambda_k)_{1\leq k\leq N}$ must be of order at least $2L_J^{(b)}$. SHT denotes spherical harmonics transform, computed from the samples of $X$ and of the $\mathcal{Y}_{\ell,m}$'s thanks to (1a). Double arrows denotes $J$ operations.

Since the needlet coefficients at a given scale $j$ depend only on a finite number of values of the function $f$, it is possible to generalize this notion to an arbitrary (possibly random) finite sequence $(f_k)_{1\leq k\leq N_j} \in \mathbb{R}^{N_j}$. The *needlet coefficients* of such sequence are the quantities

$$(\lambda_{j,k})^{-1}\sum_{p=1}^{N_j} \lambda_{j,p}\psi_{j,k}(\xi_{j,p}) f_p = \sum_{\ell\geq 0}\sum_{m=-\ell}^{\ell} b_{j,\ell}\mathcal{Y}_{\ell,m}(\xi_{j,k}) \sum_{p=1}^{N_j} \lambda_{j,p}\mathcal{Y}_{\ell,m}(\xi_{j,p}) f_p.$$

If $f_k = f(\xi_{j,k})$ for some $f \in \bigoplus_{\ell=0}^{L_j}$ and $\xi_k$ are the points of a pixelization of order at least $2L_j$, then the above expressions are equal to $\gamma_{j,k}$.

Let us give the first properties of the needlet coefficients of a random field at a given scale $j \in \mathbb{N}$. Let $X$ be a stationary field, $X_j$ like in Eq. (3) with $\mathbf{B}_{j,\ell} = 1$ if $\ell \leq L_j^{(b)}$ and $\mathbf{B}_{j,\ell} = 0$ if $\ell > 2L_j^{(b)}$, and $(\xi_{j,k}, \lambda_{j,k})_{1\leq k\leq N_j}$ a pixelization of order $4L_j^{(b)}$. The needlet coefficients of $X$ are denoted $\eta_{j,k}$. They are also the needlet coefficients of $X_j$ since $X$ and $X_j$ have the same spherical harmonics coefficients up to the frequency $\ell = L_j^{(b)}$. In the presence of an additive noise $Z_{j,k}$ in the observation of $X(\xi_{j,k})$, the "observed" needlet coefficients computed by (4) from $X + Z$ write $\eta_{j,k} + \zeta_{j,k}$, where the coefficients $\zeta_{j,k}$ are the needlet coefficients of $Z$. The next results provide the covariance structure of those coefficients at scale $j$. In our model, $X$ and $Z$ are supposed Gaussian. In this case, their needlets coefficients are Gaussian too.

**Proposition 3.** *Denote $(C_\ell)_{\ell\geq 0}$ the power spectrum of $X$. Its needlet coefficients $\eta_{j,k}$ are centered, with covariances given by*

$$\mathbf{Cov}\left[\eta_{j,k}, \eta_{j,k'}\right] = \sum_{\ell\geq 0} (b_{j,\ell})^2 C_\ell L_\ell(\xi_{j,k} \cdot \xi_{j,k'}).$$

**Proposition 4.** *Assume $Z$ of the form $Z_{j,k} = \sigma_{j,k}U_{j,k}$ where the $U_{j,k}$ are uncorrelated, centered and unit variance random variables. The needlet coefficients $\zeta_{j,k}$ of $Z$ are centered, with covariances given by*

$$\mathbf{Cov}\left[\zeta_{j,k}, \zeta_{j,k'}\right] = \frac{1}{\sqrt{\lambda_{j,k}\lambda_{j,k'}}} \sum_{p=1}^{N_j} (\lambda_{j,p}\sigma_{j,p})^2 \psi_{j,k}(\xi_{j,p})\psi_{j,k'}(\xi_{j,p}).$$

### 2.4.2. $\mathcal{B}$-adic needlets

In this paper we shall fix some constant $\mathcal{B} > 1$ and consider $\mathcal{B}$-adic window functions.

**Assumption 1.** *There exist $M \geq 3$ and a $M$-differentiable real function $\mathbf{a}$ supported in $[-\mathcal{B}, \mathcal{B}]$ and identically equal to 1 on $[-\mathcal{B}^{-1}, \mathcal{B}^{-1}]$ such that*

$$b_{j,\ell} = \mathbf{b}\left(\mathcal{B}^{-j}\ell\right)$$

*where $\mathbf{b}(\cdot) = \mathbf{a}(\cdot/\mathcal{B}) - \mathbf{a}(\cdot)$.*

For such window functions, $L_j^{(b)} = \mathcal{B}^{j+1}$. These spectral windows are not as general as those used by [5]. Indeed, it has been shown by [9] that one can take advantage of the relaxation of the $\mathcal{B}$-adic scheme originally proposed in the definition of needlets to optimize their non-asymptotic localization properties. In the following, since we are concerned with asymptotic properties, we will make use of the $\mathcal{B}$-adic structure of Assumptions 1 and 2, so that the spatial localization property of the needlet takes the convenient form of the next proposition.

**Proposition 5** ([17]). *There exists a constant $c > 0$ which depends only on the function $\mathbf{b}$ such that for all $j \in \mathbb{N}$, $k \in \{1, \ldots, N_j\}$ and $\xi \in \mathbb{S}$*

$$|\psi_{j,k}(\xi)| \leq \frac{c\mathcal{B}^j}{\left(1 + \mathcal{B}^j d(\xi, \xi_{j,k})\right)^M} \quad .$$

The stochastic counterpart of this analytical result is that the needlet coefficients of a stationary field are asymptotically uncorrelated as $j \to \infty$, except for points at a distance of order $\mathcal{B}^{-j}$ or less. For this purpose and throughout this article, we make on the power spectrum of $X$ the same following regularity assumption as in [1, 12, 11, 13].

**Assumption 2.** *There exist $\alpha > 2$ and a sequence of functions $\mathbf{g}_j : [\mathcal{B}^{-1}, \mathcal{B}] \to \mathbb{R}$, $j \in \mathbb{N}$, such that*

$$C_\ell = \ell^{-\alpha} \mathbf{g}_j \left(\mathcal{B}^{-j}\ell\right)$$

*for every $\ell \in [\mathcal{B}^{j-1}, \mathcal{B}^{j+1}]$. Moreover, there exist positive numbers $c_0, \ldots, c_M$ such that for all $j \in \mathbb{N}$, $c_0^{-1} \leq \mathbf{g}_j \leq c_0$ and for all $r \leq M$, $\sup_{\mathcal{B}^{-1} \leq u \leq \mathcal{B}} \left|\frac{d^r}{du^r} \mathbf{g}_j(u)\right| \leq c_r$.*

**Proposition 6** ([1]). *Let $X$ be stationary with a power spectrum satisfying Assumption 2 and $\eta_{j,k}$ its needlet coefficients. Then there exist a constant $c > 0$ such that for all $j \in \mathbb{N}$ and $k, k' \in \{1, \ldots, N_j\}$*

$$|\mathbf{Cor}\left[\eta_{j,k}, \eta_{j,k'}\right]| \leq \frac{c}{\left(1 + \mathcal{B}^j d(\xi_{j,k}, \xi_{j,k'})\right)^M} \quad .$$

**Remark.** *A generalization of this $\mathcal{B}$-adic framework, in a different direction to the one of [9] can be found in [7, 13, 11]. The authors do not suppose that the function $\mathbf{a}$ (or $\mathbf{b}$) is compactly supported and obtain localization and asymptotic uncorrelation results similar to Propositions 5 and 6.*

## 3. Estimation results

In this Section, we present a new procedure for the estimation of the angular power spectrum of $X$ in the model of Eq. (2) based on the needlet coefficients of $Y$. The properties of needlets described in Section 2.4 allow to take into account the local signal-to-noise ratio in the estimation of the (however) global spectrum. This spatial accuracy is at the cost of a lower frequential precision: not every value of $C_\ell$ will be estimated, but only the mean values in the bands defined by the windows $b_{j,\ell}$.

### 3.1. Parameters and estimators

We shall provide below estimators for the parameters

$$C^{(j)} \stackrel{\text{def}}{=} (4\pi)^{-1} \sum_{\ell \geq 0} (2\ell + 1) (b_{j,\ell})^2 C_\ell$$

which are smoothed versions of the power spectrum of $X$ at successive scales. A comparison of a Riemannian sum to an integral shows that, under Assumption 2, $c_1 \mathcal{B}^{j(2-\alpha)} \leq C^{(j)} \leq c_2 \mathcal{B}^{j(2-\alpha)}$ for some $c_1, c_2 > 0$.

A sequence of estimators $\widehat{C}^{(j)}$ is said *consistent* (in the high-frequency asymptotics) if $\mathcal{E}(\widehat{C}^{(j)}) \underset{j \to \infty}{\longrightarrow} 0$, where

$$\mathcal{E}(\widehat{C}^{(j)}) \stackrel{\text{def}}{=} \frac{\mathbf{E}\left(\widehat{C}^{(j)} - C^{(j)}\right)^2}{\left(C^{(j)}\right)^2} \quad .$$

The Needlet Spectral Estimators (NSE) are defined by

$$\widehat{C}^{(j)} \stackrel{\text{def}}{=} \sum_{k=1}^{N_j} w_{j,k} \left[ (\gamma_{j,k})^2 - (n_{j,k})^2 \right] \tag{5}$$

with some weights such that $\sum_{k=1}^{N_j} w_{j,k} = 1$ and where

$$n_{j,k} \stackrel{\text{def}}{=} \left( \frac{1}{\lambda_{j,k}} \sum_{p=1}^{N_j} (\lambda_{j,p} \sigma_{j,p} \psi_{j,k}(\xi_p))^2 \right)^{1/2} \tag{6}$$

is the standard deviation of the needlet coefficients of the masked noise (it follows from Proposition 4 with $k = k'$ and $W_{j,k} \sigma_{j,k}$ in place of $\sigma_{j,k}$). The $\gamma_{j,k}$'s are the observed needlet coefficients of $Y$ and $n_{j,k}$ are defined in Eq. (6). The weights should ideally be designed according to the local noise level and the local effect of the mask. For this purpose, given a sequence $(t_j)_{j \in \mathbb{N}}$ in $\mathbb{R}$, let us define the set of *kept coefficients*

$$\mathcal{K}_j \stackrel{\text{def}}{=} \left\{ k \in \{1, \ldots, N_j\} : \left\{ \sum_{p=1}^{N_j} \lambda_{j,p} (1 - W_{j,p})^2 (\psi_{j,k}(\xi_{j,p}))^2 \right\}^{1/2} \leq t_j \right\}.$$

As the second-order characteristics of the noise are known, its effect on the square of the needlet coefficients is completely subtracted in (5). Thus, the bias of the NSE is caused only by the mask.

**Proposition 7.** *For any choice of weights such that $w_{j,k} = 0$ for $k \notin \mathcal{K}_j$, there exists a constant $c$ such that*

$$\frac{\left| \mathbf{E}\, \widehat{C}^{(j)} - C^{(j)} \right|}{C^{(j)}} \leq c \mathcal{B}^{\alpha j} t_j.$$

*In particular, if $W_j \equiv 1$ (no mask), then for any choice of weights, $\widehat{C}^{(j)}$ is unbiased.*

This result is completed by the following one about the variance of $\widehat{C}$ when there is no noise.

**Proposition 8.** *Suppose that $\sigma_{j,k} = 0$ for all $k \in \{1, \ldots, N_j\}$ and that $w_{j,k} = 0$ for $k \notin \mathcal{K}_j$. Then there exists $c > 0$ such that*

$$\frac{\mathbf{Var}\left(\widehat{C}^{(j)}\right)}{C^2} \leq c \left\{ \sum_{k,k'=1}^{N_j} \frac{|w_{j,k} w_{j,k'}|}{\left(1 + \mathcal{B}^j d(\xi_{j,k}, \xi_{j,k'})\right)^{2M}} + \frac{N_j t_j}{C^{(j)}} + \left(\frac{N_j t_j}{C^{(j)}}\right)^2 \right\}$$

If the needlets coefficients at a scale $j$ were independent, centered Gaussian variables with unknown variance $C^{(j)}$ and observed with independent centered heteroscedastic Gaussian additive errors of variance $(n_{j,k})^2$, then the Maximum Likelihood of $C^{(j)}$ would be reached for the value $\widehat{C}^{(j)}$ given by (5) with weights $w_{j,k} \propto \left(C^{(j)} + (n_{j,k})^2\right)^{-2}$ (see *e.g.* [5]). In reality, the uncorrelation of the needlets coefficients holds only approximately in the sens of Proposition 6. However, if some previous, possibly rough, estimates for $C^{(j)}$ are available, say $\overline{C}^{(j)}$, the following weights are still meaningful.

$$w_{j,k} = \frac{1}{S} \left( \overline{C}^{(j)} + (n_{j,k})^2 \right)^{-2} \mathbf{1}_{k \in \mathcal{K}_j} \tag{7}$$

where $S \stackrel{\text{def}}{=} \sum_{k=1}^{N_j} \left( \overline{C}^{(j)} + (n_{j,k})^2 \right)^{-2} \mathbf{1}_{k \in \mathcal{K}_j}$.

### 3.2. Consistency of the estimators

Suppose that Assumptions 1 and 2 are satisfied. Theorems 1 and 2 below provide conditions for the consistency of the needlet spectral estimators (5). Theorem 1 deals with generic weights whereas in Theorem 2 we consider weights of the form (7).

**Theorem 1.** *Assume that*

*i) the weights $w_{j,k}$ are such that $w_{j,k} = 0$ for $k \notin \mathcal{K}_j$ and that moreover*

$$\sum_{k,k'=1}^{N_j} \frac{|w_{j,k} w_{j,k'}|}{\left(1 + \mathcal{B}^j d(\xi_{j,k}, \xi_{j,k'})\right)^M} \xrightarrow[j \to \infty]{} 0$$

*and*

$$\sup_{1 \leq p \leq N_j} \sum_{k,k'=1}^{N_j} \frac{|w_{j,k} w_{j,k'}|}{(1 + \mathcal{B}^j d(\xi_{j,p}, \xi_{j,k}))^{2M}(1 + \mathcal{B}^j d(\xi_{j,p}, \xi_{j,k'}))^{2M}} \underset{j \to \infty}{=} O\left(\mathcal{B}^{-4j}\right) \ ;$$

*ii) the noise verifies*

$$\sum_{p=1}^{N_j} \lambda_{j,p} (W_{j,p} \sigma_{j,p})^4 \underset{j \to \infty}{=} o\left(\mathcal{B}^{(4-\alpha)j}\right) \ ;$$

*iii) the thresholds $t_j$ defining the set of kept coefficients are chosen such that*

$$t_j \underset{j \to \infty}{=} o\left(\mathcal{B}^{-\alpha j}\right) \ .$$

*Then the sequence $\widehat{C}^{(j)}$ is consistent.*

**Remark 1.** *Condition i) prohibits weights which are too much concentrated on few points k. For instance, thanks to Lemma 1, Condition i) is trivially satisfied if*

$$\max_{1 \leq k \leq N_j} |w_{j,k}| \underset{j \to \infty}{=} O\left(N_j^{-1}\right).$$

**Theorem 2.** *Suppose that the weights follow Eq. (7) and that*

*i) there exists a constant $\kappa_1 > 1$ and a sequence $r_j$ in $\mathbb{R}_+$ such that*

$$\frac{(1 + r_j)^2}{\text{Card}\left\{k \in \mathcal{K}_j : (n_{j,k})^2 \leq r_j C^{(j)}\right\}} \underset{j \to \infty}{=} O(\mathcal{B}^{-\kappa_1 j}) \ ;$$

*ii) there exists a constant $\kappa_2 \in \mathbb{R}$ such that*

$$\sum_{p=1}^{N} \lambda_p (W_{j,p} \sigma_{j,p})^2 \underset{j \to \infty}{=} O(\mathcal{B}^{\kappa_2 j})$$

*iii) the thresholds $t_j$ defining the set of kept coefficients are such that*

$$t_j \underset{j \to \infty}{=} o\left(\mathcal{B}^{-\alpha j}\right) \ ;$$

*iv) there exists a constant $c > 0$ such that for all $j \in \mathbb{N}$*

$$c^{-1} C^{(j)} \leq \overline{C}^{(j)} \leq c\, C^{(j)} \ .$$

*Then, for $M > \max\left\{\frac{\kappa_2 + \alpha}{\kappa_1 - 1} \ ; \ \frac{5 + \alpha - \kappa_1 + \kappa_2}{\kappa_1 - 1}\right\}$, the sequence $\widehat{C}^{(j)}$ is consistent.*

The hypotheses of Theorem 2 can be interpreted in the following way: Condition i) describes a trade-off between an arbitrary upper bound on the noise variance and the number of observed points under this bound. The noise level can be as high as wanted (in the – very weak – limit of Condition ii)) in some regions, provided that it is under control at least in some other regions. An example is given below. Conditions ii) and iii) are not restrictive, neither is the condition on the number $M$ of derivatives of **a**, since one is free to chose a infinitely differentiable **a**. On the contrary, Condition iv) assumes a previous partial knowledge on the parameters $C^{(j)}$ to estimate, which is equivalent to the knowledge of the parameter $\alpha$.

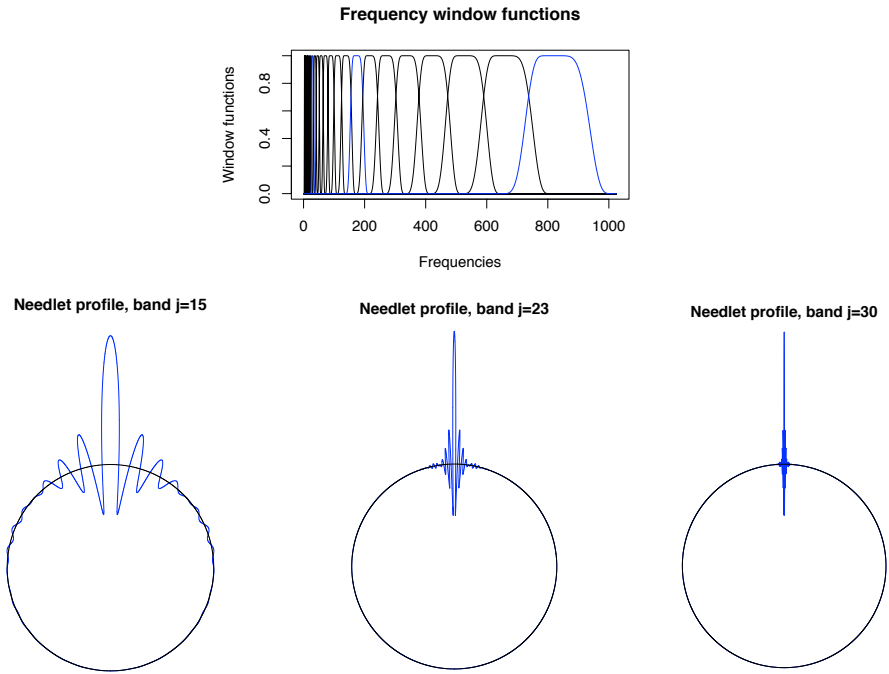

FIG 1. *Needlets in frequency and spatial spaces. Top: window functions $b_{j,\ell}$ as functions of frequency $\ell$ (in blue, the three ones corresponding to profiles plotted below). Bottom: Polar plots of the function $\theta \in [0, 2\pi[ \mapsto \sum_{\ell \geq 0} b_{j,\ell} L_\ell(\cos \theta)$ for $j = 15$, $23$ and $30$, to illustrate their localization property.*

**Example.** Suppose that there exists a $\delta > 0$ such that for each scale $j$, data with a noise variance less than $\sqrt{N_j}$ are available at least within a spherical disc of radius $\delta$. More precisely, there exist $\xi_j^* \in \mathbb{S}$ and $r_j = O\left(\left(\sqrt{N_j}\right)^\beta\right)$, $\beta < 1$ such that for all $k$ with $d(\xi_{j,k}, \xi_j^*) \leq \delta$

- $W(\xi_{j,k}) = 1$
- and $\sigma_{j,k} \leq r_j C^{(j)}$.

Then condition i) is satisfied.

## 4. Simulations

In this section, we investigate the non-asymptotic numerical performances of the NSE, in a relatively realistic setting. This is for illustration only, as there are many free parameters in the models (mask, noise level) that are only asymptotically constrained by the previous theory.

### 4.1. Settings

The constant determining the width of the bands is set to $\mathcal{B} = 1.25$. The function **a** of Assumption 1 is designed on $[\mathcal{B}^{-1}, \mathcal{B}]$ by a polynomial of degree 19 connecting 1 to 0 with $M = 9$ continuous derivatives, i.e. $\mathbf{a}(\mathcal{B}^{-1}) = 1$, $\mathbf{a}(\mathcal{B}) = 0$ and $\mathbf{a}^{(r)}(\mathcal{B}^{\pm 1}) = 0$, $1 \leq r \leq M$. Frequency window functions $b_{j,\ell}$ and some of the associated needlets are displayed in Fig. 1.

In the analysis, the first 10 scales ($0 \leq j \leq 9$) are not considered since the support of the corresponding window functions is of width less or equal than 2. These scales correspond to frequencies $\ell \leq 8$. The range of studied scales is $10 \leq j \leq 30$.

We use the HEALPix pixelization and associated direct and inverse spherical harmonic transform [8]. Needlet coefficients are computed from three experiments providing maps at respective

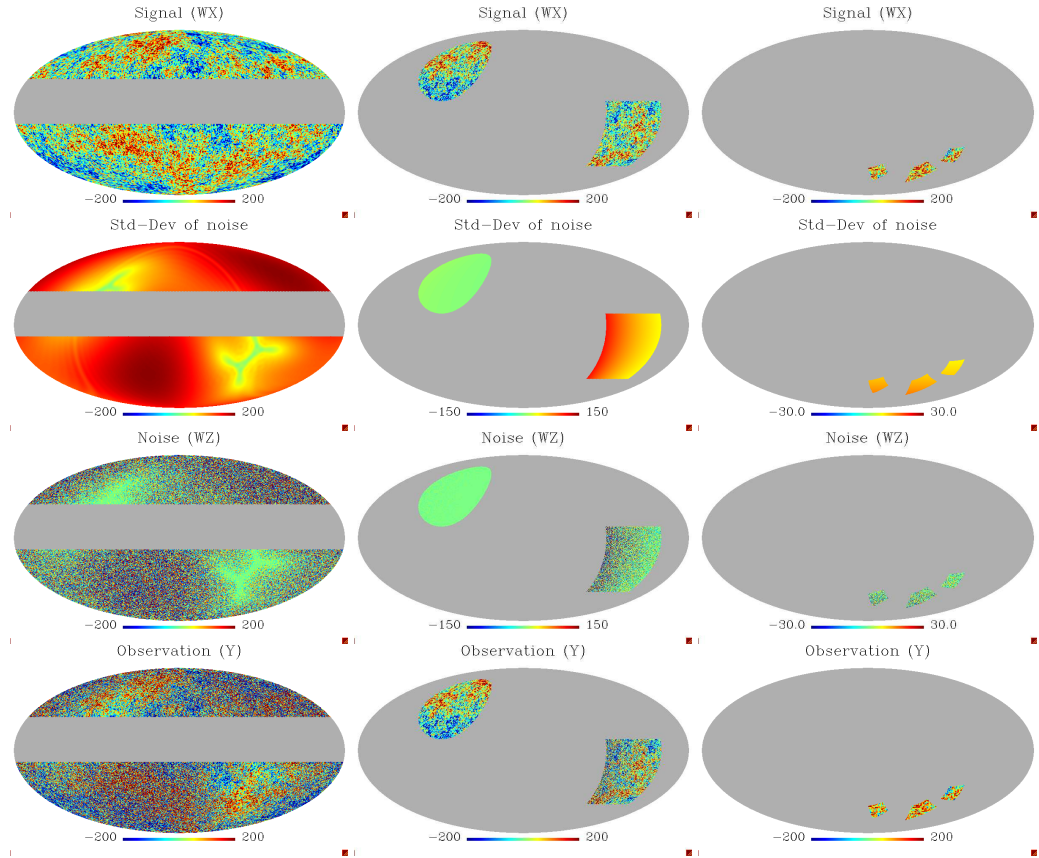

FIG 2. *Simulated data. Fist line: $WX$. Second line: $W\sigma$. Third line: $WZ$. Fourth line: $Y$. The columns correspond to the settings for the three experiments: A (fist column) used for $10 \leq j \leq 23$; B (second column) used for $23 < j \leq 26$; C (third column) used for $26 < j \leq 30$.*

HEALPix resolutions nside $= 128, 256$ and $512$ (*i.e.* the number of cubature points $N$ is 196.608, 786.432 and 3.145.728).

The simulated data are illustrated in Fig. 2. A first "experiment" (hereafter, $A$), *i.e.* a mask an a noise level map, is used for scales $10 \leq j \leq 23$. The noise level is reminiscent to the expected level form the forthcoming Planck experiment of the European Spatial Agency[1]. Two other experiments (hereafter, $B$ and $C$), with synthetic masks and noise-level maps, are used respectively for scales $23 < j \leq 26$ and $26 < j \leq 30$.

### 4.2. Results

The distribution on the NSE estimators for first the scales $10 \leq j \leq 30$ is estimated by 800 Monte Carlo replicates in the setting of previous section. Figure 3 illustrates the convergence of the sequence $\widehat{C}^{(j)}$. The approximate normal distribution of the estimator is highlighted by the quantile-to-quantile plots of Figure 4 and the $p$-values from Anderson-Darling goodness-of-fit test.

## 5. Conclusion

In this work, we have proved the mean-square consistency of a spherical spectral estimator based on the needlet analysis of a stationary field with missing data and non-homogeneous additive noise. Our model is relatively realistic and could be enriched to fit the requirement of particular applications. A possible development of this work is the study of the rate of convergence, limiting law, efficiency, etc. of those estimators.

---
[1]see www.rssd.esa.int/Planck/

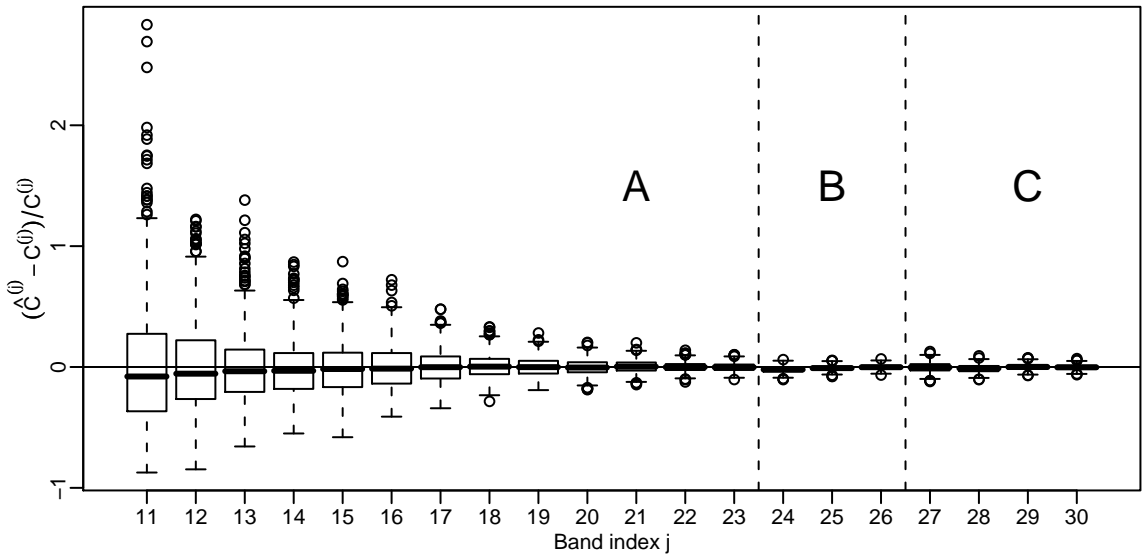

Fig 3. *Box-and-whisker plots of the NSE estimator $\widehat{C}^{(j)}$. Regions A,B, and C delimits the range of the three experiments, from the coarsest (and wide coverage) to the finest (and small coverage) one.*

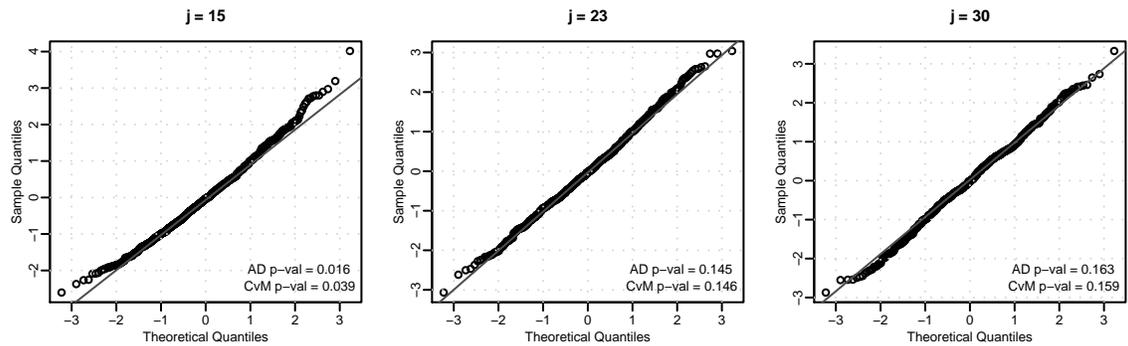

Fig 4. *Distribution of the standardized NSE estimator $\widehat{C}^{(j)}$ for bands 15, 23 and 30. The p-values are obtained from the Anderson-Darling and Cramér-von Mises tests of Gaussianity. At low frequencies, the estimator is significatively skewed.*

Although this has not be stated here for brevity, it must be noticed that this method allows easy adaptation to a context with longitudinal data, *i.e.* of the form

$$Y_{j,k,e} = W_{j,k,e}\left[X\left(\xi_{j,k}\right) + Z_{j,k,e}\right], \ j \in \mathbb{N}, \ k \in \{1,\ldots,N_j\}, \ e \in \{1,\ldots,E_j\}.$$

for some sequence of integer $\{E_j\}_{j\in\mathbb{N}}$. This generalization of model (2) is of high importance for practical applications, when for each scale $j$ multiple independent observations (indexed by $e$) of the signal of interest are available, such those provided by the various CMB experiments (see [5] for details).

### *Acknowledgments*

We are grateful to P. Baldi, J.-F. Cardoso, J. Delabrouille, G. Kerkyacharian, D. Marinucci and and D. Picard for usefull discussions. The simulations of Section 4 have been conducted using the HEALPix package [8].

## 6. Proofs

In this section, for sake of notational simplicity,

- the sub- and superscripts $j$ are omitted
- for quantities $a(j)$ and $b(j)$ depending on $j$, $a \lesssim b$ means: $\exists c > 0, \forall j \in \mathbb{N}, a(j) \leq cb(j)$
- $a \asymp b$ means: $a \lesssim b$ and $b \lesssim a$.

The following fact will be used extensively in the proofs. Let $(\xi_p, \lambda_p)_{1 \leq p \leq N}$ be a pixelization of order $4\mathcal{B}^{j+1}$

$$\sum_{p=1}^{N} \lambda_p \psi_k^2(\xi_p) = \int_{\mathbb{S}} \psi_k^2(\xi) d\xi = \lambda_k^{-1} \sum_{\ell,\ell' \geq 0} \mathbf{b}(\mathcal{B}^{-j}\ell)\mathbf{b}(\mathcal{B}^{-j}\ell') \int_{\mathbb{S}} L_\ell(\xi_k \cdot \xi) L_{\ell'}(\xi_k \cdot \xi) d\xi$$

$$= \lambda_k^{-1} \sum_{\ell \geq 0} \mathbf{b}^2(\mathcal{B}^{-j}\ell) \frac{2\ell+1}{4\pi} \asymp \lambda_k^{-1} \mathcal{B}^j \sum_{\ell \geq 0} \mathbf{b}^2(\mathcal{B}^{-j}\ell)\mathcal{B}^{-j}\ell \asymp 1 \qquad (8)$$

as $\int_{\mathbb{R}_+} \mathbf{b}(x)x dx > 0$ and $\lambda_k^{-1} \mathcal{B}^j \asymp \mathcal{B}^{-j}$ using (1b)-(1c).

### *6.1. Proof of Propositions 3 and 4*

$$\begin{aligned}
\mathbf{E}\left[\eta_k \eta_{k'}\right] &= (\lambda_k \lambda_{k'})^{-1/2} \mathbf{E} \iint_{\mathbb{S}} X(\xi)X(\xi')\psi_k(\xi)\psi_{k'}(\xi')d\xi d\xi' \\
&= (\lambda_k \lambda_{k'})^{-1/2} \mathbf{E} \iint_{\mathbb{S}} \sum_{\ell,\ell' \geq 0'} \sum_{m=-\ell}^{\ell} \sum_{m'=-\ell'}^{\ell'} a_{\ell,m} a_{\ell',m'} \mathcal{Y}_{\ell,m}(\xi) \mathcal{Y}_{\ell',m'}(\xi')\psi_k(\xi)\psi_{k'}(\xi')d\xi d\xi' \\
&= (\lambda_k \lambda_{k'})^{-1/2} \sum_{\ell \geq 0} C_\ell \sum_{m=-\ell}^{\ell} \int_{\mathbb{S}} \mathcal{Y}_{\ell,m}(\xi)\psi_k(\xi)d\xi \int_{\mathbb{S}} \mathcal{Y}_{\ell,m}(\xi)\psi_{k'}(\xi)d\xi \\
&= \sum_{\ell \geq 0} C_\ell b_\ell^2 \sum_{m=-\ell}^{\ell} \mathcal{Y}_{\ell,m}(\xi_k)\mathcal{Y}_{\ell,m}(\xi_{k'})
\end{aligned}$$

which proves Proposition 3. On the other hand $\zeta_k = \sum_{\ell \geq 0} \sum_{m=-\ell}^{\ell} b_\ell z_{\ell,m} \mathcal{Y}_{\ell,m}(\xi_k)$ where

$$z_{\ell,m} \stackrel{\text{def}}{=} \sum_{p=1}^{N} \lambda_p Z_p \mathcal{Y}_{\ell,m}(\xi_p)$$

Thus

$$\begin{aligned}
\mathbf{Cov}\left[\zeta_k \zeta_{k'}\right] &= \sum_{\ell,\ell' \geq 0} \sum_{m=-\ell}^{\ell} \sum_{m'=-\ell'}^{\ell'} b_\ell b_{\ell'} \mathbf{Cov}\left[z_{\ell,m}, z_{\ell',m'}\right] \mathcal{Y}_{\ell,m}(\xi_k) \mathcal{Y}_{\ell',m'}(\xi_{k'}) \\
&= \sum_{\ell,\ell' \geq 0} \sum_{m=-\ell}^{\ell} \sum_{m'=-\ell'}^{\ell'} b_\ell b_{\ell'} \left(\sum_{p=1}^{N} \lambda_p^2 \sigma_p^2 \mathcal{Y}_{\ell,m}(\xi_p) \mathcal{Y}_{\ell',m'}(\xi_p)\right) \mathcal{Y}_{\ell,m}(\xi_k) \mathcal{Y}_{\ell',m'}(\xi_{k'}) \\
&= \frac{1}{\sqrt{\lambda_k \lambda_{k'}}} \sum_{p=1}^{N} \lambda_p^2 \sigma_p^2 \psi_k(\xi_p) \psi_{k'}(\xi_p)
\end{aligned}$$

### 6.2. Proof of Proposition 7

$$\begin{aligned}
\frac{|\mathbf{E}\widehat{C} - C|}{C} &= \frac{1}{C} \left| \sum_{k \in \mathcal{K}} w_k \left\{ \mathbf{Var}\left(\gamma_k\right) - n_k^2 - C \right\} \right| \\
&= \frac{1}{C} \left| \sum_{k \in \mathcal{K}} w_k \left\{ \mathbf{Var}\left(\eta_k\right) - C \right\} \right| \\
&\leq \frac{1}{C} \sum_{k \in \mathcal{K}} w_k \lambda_k^{-1} \left| \mathbf{E}\left(\sum_{p=1}^{N} \lambda_p W_p X(\xi_p) \psi_k(\xi_p)\right)^2 - \mathbf{E}\left(\sum_{p=1}^{N} \lambda_p X(\xi_p) \psi_k(\xi_p)\right)^2 \right| \\
&= \frac{1}{C} \sum_{k \in \mathcal{K}} w_k \lambda_k^{-1} \left| \mathbf{E}\left(\sum_{p=1}^{N} \lambda_p (1 - W_p) X(\xi_p) \psi_k(\xi_p) \sum_{p=1}^{N} \lambda_p (1 + W_p) X(\xi_p) \psi_k(\xi_p)\right) \right| \\
&\leq \frac{2}{C} \sum_{k \in \mathcal{K}} w_k \lambda_k^{-1} \left\{ \mathbf{E}\left(\sum_{p=1}^{N} \lambda_p (1 - W_p) X(\xi_p) \psi_k(\xi_p)\right)^2 \mathbf{E}\left(\sum_{p=1}^{N} \lambda_p X(\xi_p) \psi_k(\xi_p)\right)^2 \right\}^{1/2} \\
&\leq \frac{8\pi \mathbf{Var}(X)}{C} \sum_{k \in \mathcal{K}} w_k \lambda_k^{-1} \left\{ \sum_{p=1}^{N} \lambda_p (1 - W_p)^2 \psi_k^2(\xi_p) \sum_{p=1}^{N} \lambda_p \psi_k^2(\xi_p) \right\}^{1/2} \\
&\lesssim \frac{N}{C} t \sum_{k \in \mathcal{K}} w_k \\
&\asymp \mathcal{B}^{\alpha j} t \quad .
\end{aligned}$$

In particular, if $W \equiv 1$, one can take $t = 0$ in the definition of $\mathcal{K}$ and the above result shows that $\widehat{C}$ is unbiased.

### 6.3. Proof of Proposition 8

Define $\tilde{\eta}_k \stackrel{\text{def}}{=} \lambda_k^{-1/2} \langle \psi_k, X \rangle_{\mathbb{H}}$ .

$$\begin{aligned}
\frac{\mathbf{Var}(\widehat{C})}{C^2} &= \frac{1}{C^2}\sum_{k,k'\in\mathcal{K}} w_k w_{k'} \mathbf{Cov}\,[\eta_k^2,\eta_{k'}^2] = \frac{2}{C^2}\sum_{k,k'\in\mathcal{K}} w_k w_{k'} \mathbf{Cov}^{\,2}[\eta_k,\eta_{k'}] \\
&= \frac{2}{C^2}\sum_{k,k'\in\mathcal{K}} w_k w_{k'} \big\{\mathbf{Cov}\,[\tilde\eta_k,\tilde\eta_{k'}] - \mathbf{Cov}\,[\tilde\eta_k,\tilde\eta_{k'}-\eta_{k'}] \\
&\qquad - \mathbf{Cov}\,[\tilde\eta_k-\eta_k,\tilde\eta_{k'}] + \mathbf{Cov}\,[\tilde\eta_k-\eta_k,\tilde\eta_{k'}-\eta_{k'}]\big\}^2 \\
&\asymp \frac{1}{C^2}\sum_{k,k'\in\mathcal{K}} w_k w_{k'}\big\{\mathbf{Cov}^{\,2}[\tilde\eta_k,\tilde\eta_{k'}] + \mathbf{Cov}^{\,2}[\tilde\eta_k,\tilde\eta_{k'}-\eta_{k'}] + \mathbf{Cov}^{\,2}[\tilde\eta_k-\eta_k,\tilde\eta_{k'}-\eta_{k'}]\big\} \quad . \quad (9)
\end{aligned}$$

For the first term in the r.h.s. of Eq. (9), we use the following bound, proved in [1].

$$\sum_{\ell\geq 0} b_\ell^2 C_\ell L_\ell(\xi_k\cdot\xi_{k'}) \lesssim \frac{\mathcal{B}^{j(2-\alpha)}}{(1+\mathcal{B}^j d(\xi_k,\xi_{k'}))^M} \quad . \qquad (10)$$

Thus, with Lemma 1 and using that $C \asymp \mathcal{B}^{j(2-\alpha)}$ and $N \asymp \mathcal{B}^{2j}$

$$\begin{aligned}
\sum_{k,k'\in\mathcal{K}} w_k w_{k'} \mathbf{Cov}^{\,2}[\tilde\eta_k,\tilde\eta_{k'}] &= \sum_{k,k'\in\mathcal{K}} w_k w_{k'} \Big\{\sum_{\ell\geq 0} b_\ell^2 C_\ell L_\ell(\xi_k\cdot\xi_{k'})\Big\}^2 \\
&\lesssim \mathcal{B}^{2j(2-\alpha)} \sum_{k,k'\in\mathcal{K}} \frac{|w_k w_{k'}|}{(1+\mathcal{B}^j d(\xi_k,\xi_{k'}))^{2M}} \\
&\lesssim C^2 \sum_{k,k'\in\mathcal{K}} \frac{|w_k w_{k'}|}{(1+\mathcal{B}^j d(\xi_k,\xi_{k'}))^{2M}} \quad .
\end{aligned}$$

For the other terms of the decomposition (9), one can write

$$\begin{aligned}
\mathbf{Var}(\tilde\eta_k-\eta_k) &= \lambda_k^{-1}\mathbf{E}\Big(\sum_{p=1}^N \lambda_p(1-W_p)X(\xi_p)\psi_k(\xi_p)\Big)^2 \\
&\leq \lambda_k^{-1}\mathbf{E}\Big(\sum_{p=1}^N \lambda_p X^2(\xi_p) \sum_{p=1}^N \lambda_p(1-W_p)\psi_k^2(\xi_p)\Big) \\
&= \lambda_k^{-1} 4\pi \mathbf{Var}(X) \sum_{p=1}^N \lambda_p(1-W_p)\psi_k^2(\xi_p) \\
&\lesssim Nt \qquad\qquad (11)
\end{aligned}$$

thus

$$\begin{aligned}
\sum_{k,k'\in\mathcal{K}} w_k w_{k'}\{\mathbf{Cov}\,[\tilde\eta_k-\eta_k,\tilde\eta_{k'}]\}^2 &\leq \sum_{k\in\mathcal{K}} w_k \mathbf{Var}(\tilde\eta_k-\eta_k) \sum_{k'\in\mathcal{K}} w_{k'} \mathbf{Var}(\tilde\eta_k) \\
&\lesssim CNt
\end{aligned}$$

and

$$\begin{aligned}
\sum_{k,k'\in\mathcal{K}} w_k w_{k'}\{\mathbf{Cov}\,[\tilde\eta_k-\eta_k,\tilde\eta_{k'}-\eta_{k'}]\}^2 &\leq \Big\{\sum_{k\in\mathcal{K}} w_k \mathbf{Var}(\tilde\eta_k-\eta_k)\Big\}^2 \\
&\lesssim N^2(t)^2 \quad .
\end{aligned}$$

Putting these terms together establishes Proposition 8.

### 6.4. Proof of Theorems 1 and 2

We have
$$\mathcal{E}(\widehat{C}) = \frac{|\mathbf{E}\,\widehat{C} - C|^2}{C^2} + \frac{\mathbf{Var}\,(\widehat{C})}{C^2} \quad .$$

By Proposition 7, the bias term is bounded by
$$\frac{|\mathbf{E}\,\widehat{C} - C|}{C} \lesssim \mathcal{B}^{\alpha j} t$$

and converges to 0 thanks to Condition iii).

On the other hand
$$\begin{aligned}
\frac{\mathbf{Var}\,(\widehat{C})}{C^2} &= \frac{1}{C^2} \sum_{k,k'=1}^N w_k w_{k'} \mathbf{Cov}\,[\gamma_k^2, \gamma_{k'}^2] = \frac{2}{C^2} \sum_{k,k'=1}^N w_k w_{k'} \mathbf{Cov}\,^2[\gamma_k, \gamma_{k'}] \\
&= \frac{2}{C^2} \sum_{k,k' \in \mathcal{K}} w_k w_{k'} \left\{ \mathbf{Cov}\,[\eta_k, \eta_{k'}] + \mathbf{Cov}\,[\zeta_k, \zeta_{k'}] \right\}^2 \\
&\asymp \frac{1}{C^2} \sum_{k,k' \in \mathcal{K}} w_k w_{k'} \left\{ \mathbf{Cov}\,^2[\eta_k, \eta_{k'}] + \mathbf{Cov}\,^2[\zeta_k, \zeta_{k'}] \right\} \quad . \quad (12)
\end{aligned}$$

From Proposition 8
$$\frac{1}{C^2} \sum_{k,k' \in \mathcal{K}} w_k w_{k'} \mathbf{Cov}\,^2[\eta_k, \eta_{k'}] \lesssim \sum_{k,k'=1}^N \frac{|w_k w_{k'}|}{(1 + \mathcal{B}^j d(\xi_k, \xi_{k'}))^{2M}} + \frac{Nt}{C} + \left(\frac{Nt}{C}\right)^2$$

where $\frac{Nt}{C} \lesssim \mathcal{B}^{\alpha j} t \xrightarrow[j \to \infty]{} 0$ with Condition iii). To establish Theorems 1 and 2 it suffices now to prove that in both cases
$$\sum_{k,k'=1}^N \frac{|w_k w_{k'}|}{(1 + \mathcal{B}^j d(\xi_k, \xi_{k'}))^{2M}} \xrightarrow[j \to \infty]{} 0 \quad (13)$$

and
$$\frac{1}{C^2} \sum_{k,k' \in \mathcal{K}} w_k w_{k'} \mathbf{Cov}\,^2[\zeta_k, \zeta_{k'}] \xrightarrow[j \to \infty]{} 0 \quad . \quad (14)$$

Let us consider separately the cases of Theorems 1 and 2.

### 6.5. End of proof of Theorem 1

The convergence of Eq. (13) is nothing else but Condition i). To prove (14), use the Property (1c) and the Cauchy-Schwartz inequality and write

$$
\begin{aligned}
\sum_{k,k' \in \mathcal{K}} |w_k w_{k'}| \mathbf{Cov}^2[\zeta_k, \zeta_{k'}] &= \sum_{k,k' \in \mathcal{K}} \frac{|w_k w_{k'}|}{\lambda_k \lambda_{k'}} \Big\{ \sum_{p=1}^N \lambda_p^2 W_p^2 \sigma_p^2 \psi_k(\xi_p) \psi_{k'}(\xi_p) \Big\}^2 \\
&\asymp \sum_{k,k' \in \mathcal{K}} |w_k w_{k'}| \Big\{ \sum_{p=1}^N \lambda_p W_p^2 \sigma_p^2 \psi_k(\xi_p) \psi_{k'}(\xi_p) \Big\}^2 \\
&\lesssim \sum_{k,k' \in \mathcal{K}} |w_k w_{k'}| \Big\{ \sum_{p=1}^N \lambda_p W_p^4 \sigma_p^4 \Big\} \Big\{ \sum_{p=1}^N \lambda_p |\psi_k(\xi_p) \psi_{k'}(\xi_p)|^2 \Big\} \\
&\lesssim \Big\{ \sum_{p=1}^N \lambda_p W_p^4 \sigma_p^4 \Big\} \sum_{k,k' \in \mathcal{K}} |w_k w_{k'}| \\
&\quad \times \mathcal{B}^{2j} \sum_{p=1}^N \frac{1}{(1+\mathcal{B}^j d(\xi_p, \xi_k))^{2M}(1+\mathcal{B}^j d(\xi_p, \xi_{k'}))^{2M}} \\
&\underset{j \to \infty}{=} o(C)
\end{aligned}
$$

with Conditions i) and ii) and Proposition 5.

### 6.6. End of proof of Theorem 2

The two following remarks will prove useful. First, from Condition iv), $(c^{-1} \wedge 1)(C + n_k^2) \leq (\overline{C} + n_k^2) \leq (c \vee 1)(C + n_k^2)$ so that

$$w_k \asymp \frac{1}{S} (C + n_k^2)^{-2} \mathbf{1}_{k \in \mathcal{K}} \quad . \tag{15}$$

Second, using again Condition iv)

$$S = \sum_{k \in \mathcal{K}} (\overline{C} + n_k^2)^{-2} \geq \sum_{k \in V} (\overline{C} + n_k^2)^{-2} \geq \frac{\operatorname{Card} V}{(1+r)^2 \overline{C}^2} \gtrsim \frac{\operatorname{Card} V}{(1+r)^2 C^2} \tag{16}$$

where

$$V \stackrel{\text{def}}{=} \{k \in \mathcal{K} : n_k^2 \leq rC\} \quad . \tag{17}$$

The convergence of Eq. (13) is established thanks to Conditions i) and iii) and Lemma 1

$$
\begin{aligned}
\sum_{k,k'=1}^N \frac{|w_k w_{k'}|}{(1+\mathcal{B}^j d(\xi_k, \xi_{k'}))^{2M}} &\asymp \frac{1}{S^2} \sum_{k,k' \in \mathcal{K}} \frac{(C+n_k^2)^{-2}(C+n_{k'}^2)^{-2}}{(1+\mathcal{B}^j d(\xi_k, \xi_{k'}))^{2M}} \\
&\leq \frac{C^{-4}}{S^2} \sum_{k,k' \in \mathcal{K}} \frac{1}{(1+\mathcal{B}^j d(\xi_k, \xi_{k'}))^{2M}} \\
&\lesssim \frac{(1+r)^4}{(\operatorname{Card} V)^2} (\operatorname{Card} \mathcal{K}) \underset{j \to \infty}{\longrightarrow} 0
\end{aligned}
$$

For the convergence of Eq. (14), note that $\mathbf{Cov}[\zeta_k, \zeta_{k'}] \asymp \sum_{p=1}^N \lambda_p W_p^2 \sigma_p^2 \psi_k(\xi_p) \psi_{k'}(\xi_p)$ where the implicit constant does not depend on $k, k'$ (Proposition 4 and property (1c) of the pixelization).

For any constant $\kappa \in (0,1)$, the value of which will be fixed later, we have

$$\frac{1}{C^2}\sum_{k,k'\in\mathcal{K}} w_k w_{k'} \mathbf{Cov}^2[\zeta_k,\zeta_{k'}] \asymp \frac{1}{C^2 S^2}\sum_{k,k'\in\mathcal{K}} (C+n_k^2)^{-2}(C+n_{k'}^2)^{-2}$$

$$\times \Big\{\sum_{p=1}^N \lambda_p W_p^2 \sigma_p^2 \psi_k(\xi_p)\psi_{k'}(\xi_p)\Big\}^2$$

$$\lesssim \frac{1}{C^2 S^2}(A_1+A_2+A_3+A_4)$$

with

$$A_1 = \sum_{\substack{k,k'\in\mathcal{K}\\d(\xi_k,\xi_{k'})\leq \mathcal{B}^{-\kappa j}}} (C+n_k^2)^{-2}(C+n_{k'}^2)^{-2}\Big\{\sum_{p=1}^N \lambda_p W_p^2 \sigma_p^2 \psi_k(\xi_p)\psi_{k'}(\xi_p)\Big\}^2$$

$$A_2 = \sum_{\substack{k,k'\in\mathcal{K}\\d(\xi_k,\xi_{k'})> \mathcal{B}^{-\kappa j}}} (C+n_k^2)^{-2}(C+n_{k'}^2)^{-2}\Big\{\sum_{p\in D(k,k')} \lambda_p W_p^2 \sigma_p^2 \psi_k(\xi_p)\psi_{k'}(\xi_p)\Big\}^2$$

where $D(k,k') \stackrel{\text{def}}{=} \{p: d(\xi_k,\xi_p) > \tfrac{1}{2}\mathcal{B}^{-\kappa j} \text{ and } d(\xi_{k'},\xi_p) > \tfrac{1}{2}\mathcal{B}^{-\kappa j}\}$

$$A_3 = \sum_{\substack{k,k'\in\mathcal{K}\\d(\xi_k,\xi_{k'})> \mathcal{B}^{-\kappa j}}} (C+n_k^2)^{-2}(C+n_{k'}^2)^{-2}\Big\{n_k^2 \sum_{p:d(\xi_k,\xi_p)\leq \tfrac{1}{2}\mathcal{B}^{-\kappa j}} \lambda_p \psi_k(\xi_p)\psi_{k'}(\xi_p)\Big\}^2$$

$$A_4 = \sum_{\substack{k,k'\in\mathcal{K}\\d(\xi_k,\xi_{k'})> \mathcal{B}^{-\kappa j}}} (C+n_k^2)^{-2}(C+n_{k'}^2)^{-2}\Big\{\sum_{p:d(\xi_k,\xi_p)\leq \tfrac{1}{2}\mathcal{B}^{-\kappa j}} \lambda_p |W_p^2 \sigma_p^2 - n_k^2|\psi_k(\xi_p)\psi_{k'}(\xi_p)\Big\}^2 \; .$$

**Remark.** *As noticed in Section 3.1, the weights proposed here are linked with the asymptotic uncorrelation of the needlet coefficients. $A_1$ includes the near-diagonal terms of the noise covariance in the needlet coefficients domain ; $A_2$, $A_3$ and $A_4$ include the off-diagonal covariances, for which the localization of the needlets is crucial. In $A_2$, the "near-disjointness" of the supports of $\psi_k$ and $\psi_{k'}$ is only considered sufficiently far away from $\xi_k$ and $\xi_{k'}$, while in $A_3$ it is considered around $\xi_k$ (or equivalently $\xi_{k'}$). Finally, $A_4$ depends on the regularity of $p \mapsto \sigma_p^2$.*

We shall use the fact that: $\forall x, C \in \mathbb{R}_+$, $\frac{x^2}{(C+x^2)^2} \leq \min\left\{\frac{x^2}{C^2}\ \frac{1}{4C}\ ;\ \frac{1}{x^2}\right\}$.

With Cauchy-Schwartz inequality, $W \leq 1$ and the properties of the cubature points (including Lemma 2) we have

$$A_1 \leq \sum_{\substack{k,k'\in\mathcal{K}\\d(\xi_k,\xi_{k'})\leq \mathcal{B}^{-\kappa j}}} (C+n_k^2)^{-2}(C+n_{k'}^2)^{-2}\Big\{\sum_{p=1}^N \lambda_p \sigma_p^2 \psi_k^2(\xi_p)\Big\}\Big\{\sum_{p=1}^N \lambda_p \sigma_p^2 \psi_{k'}^2(\xi_p)\Big\}$$

$$= \sum_{k\in\mathcal{K}} \frac{n_k^2}{(C+n_k^2)^2} \sum_{\substack{k'\in\mathcal{K}\\d(\xi_k,\xi_{k'})\leq \mathcal{B}^{-\kappa j}}} \frac{n_{k'}^2}{(C+n_{k'}^2)^2}$$

$$\lesssim \frac{N^2}{C^2}\mathcal{B}^{-2\kappa j} \; .$$

Thanks to Eq. (16) and Condition i)

$$\frac{A_1}{C^2 S^2} \lesssim N^2 \frac{(1+r)^4}{(\text{Card }V)^2}\mathcal{B}^{-2\kappa j} \asymp \mathcal{B}^{2j(2-\kappa_1-\kappa)} \xrightarrow[j\to\infty]{} 0$$

provided that

$$\kappa > 2 - \kappa_1 \; . \tag{18a}$$

To give upper bounds on $A_2$ and $A_3$, we use Proposition 5. Define

$$D(k,k') \stackrel{\text{def}}{=} \left\{ p : d(\xi_k, \xi_p) > \frac{1}{2}\mathcal{B}^{-\kappa j} \text{ and } d(\xi_{k'}, \xi_p) > \frac{1}{2}\mathcal{B}^{-\kappa j} \right\}.$$

Write

$$\begin{aligned}
A_2 &\leq \sum_{k,k' \in \mathcal{K}} (C + n_k^2)^{-2}(C + n_{k'}^2)^{-2} \\
&\quad \times \Big\{ \sum_{p \in D(k,k')} \lambda_p W_p^2 \sigma_p^2 \psi_k(\xi_p) \psi_{k'}(\xi_p) \Big\}^2 \\
&\lesssim \sum_{k,k' \in \mathcal{K}} (C + n_k^2)^{-2}(C + n_{k'}^2)^{-2} \\
&\quad \times \Big\{ \mathcal{B}^{2j} \sum_{p \in D(k,k')} \frac{\lambda_p W_p^2 \sigma_p^2}{(1 + \mathcal{B}^j d(\xi_p, \xi_k))^M (1 + \mathcal{B}^j d(\xi_p, \xi_{k'}))^M} \Big\}^2 \\
&\leq \sum_{k,k' \in \mathcal{K}} (C + n_k^2)^{-2}(C + n_{k'}^2)^{-2} \Big\{ \frac{\mathcal{B}^{2j}}{(1 + \frac{1}{2}\mathcal{B}^{(1-\kappa)j})^M} \sum_{p=1}^N \lambda_p W_p^2 \sigma_p^2 \Big\}^2 \\
&\leq S^2 \left( \sum_{p=1}^N \lambda_p W_p^2 \sigma_p^2 \right)^2 \mathcal{B}^{2j(2-(1-\kappa)M)}
\end{aligned}$$

so that, with Eq. (16) and Condition ii),

$$\frac{A_2}{C^2 S^2} \lesssim \left( \sum_{p=1}^N \lambda_p \sigma_p^2 \right)^2 \mathcal{B}^{2j(\alpha-(1-\kappa)M)} \asymp \mathcal{B}^{2j(\kappa_2 + \alpha - (1-\kappa)M)} \xrightarrow[j \to \infty]{} 0$$

provided that

$$M > \frac{\kappa_2 + \alpha}{1 - \kappa} \quad . \tag{18b}$$

Remark now that if $d(\xi_k, \xi_{k'}) > \mathcal{B}^{-\kappa j}$ and $d(\xi, \xi_k) \leq \frac{1}{2}\mathcal{B}^{-\kappa j}$ then $d(\xi, \xi_{k'}) \geq \frac{1}{2}d(\xi_k, \xi_{k'})$ so that $\frac{1}{(1 + \mathcal{B}^j d(\xi, \xi_k))(1 + \mathcal{B}^j d(\xi, \xi_{k'}))} \gtrsim \frac{1}{1 + \mathcal{B}^j d(\xi_k, \xi_{k'})}$. Thus

$$\begin{aligned}
A_3 &\lesssim \sum_{\substack{k,k' \in \mathcal{K} \\ d(\xi_k, \xi_{k'}) > \mathcal{B}^{-\kappa j}}} \frac{n_k^4}{(C + n_k^2)^2 (C + n_{k'}^2)^2} \sum_{p:d(\xi_k, \xi_p) \leq \frac{1}{2}\mathcal{B}^{-\kappa j}} \lambda_p |\psi_k(\xi_p)\psi_{k'}(\xi_p)|^2 \\
&\lesssim \sum_{\substack{k,k' \in \mathcal{K} \\ d(\xi_k, \xi_{k'}) > \mathcal{B}^{-\kappa j}}} \frac{n_k^4}{(C + n_k^2)^2 (C + n_{k'}^2)^2} \mathcal{B}^{4j} \sum_{p:d(\xi_k, \xi_p) \leq \frac{1}{2}\mathcal{B}^{-\kappa j}} \frac{\lambda_p}{(1 + \mathcal{B}^j d(\xi_k, \xi_{k'}))^{2M}} \\
&\lesssim \sum_{k,k' \in \mathcal{K}} \frac{n_k^4}{(C + n_k^2)^2 (C + n_{k'}^2)^2} \mathcal{B}^{4j} \mathcal{B}^{-j(1-\kappa)2M} \mathcal{B}^{-2\kappa j} \\
&\leq \frac{N^2}{C^2} \mathcal{B}^{2j(2-(1-\kappa)M)}
\end{aligned}$$

and with Condition i),

$$\frac{A_3}{C^2 S^2} \lesssim \frac{(1+r)^4}{(\operatorname{Card} V)^2} N^2 \mathcal{B}^{2j(2-(1-\kappa)M)} \asymp \mathcal{B}^{2j(4-\kappa_1-(1-\kappa)M)} \xrightarrow[j \to \infty]{} 0$$

provided that

$$M > \frac{4 - \kappa_1}{1 - \kappa} \quad . \tag{18c}$$

For the last term $A_4$, we use Condition ii), the same arguments as for $A_3$, and Lemma 9 which takes advantage of the localization of $\psi_k$ around point $\xi_k$ to compare the local noise level in direct and needlet domains. The estimate depends on the regularity $\varrho$ of $\sigma^2$.

**Lemma 9.** *Define* $\varrho_j \stackrel{def}{=} \max_{1 \leq k < k' \leq N_j} \frac{|(W_{j,k}\sigma_{j,k})^2 - (W_{j,k'}\sigma_{j,k'})^2|}{d(\xi_{j,k}, \xi_{j,k'})}$. *There exist a constant* $c > 0$ *such that for all* $j \in \mathbb{N}$, $\epsilon > 0$, *and* $k, k' \in \{1, \ldots, N_j\}$ *with* $d(\xi_k, \xi_{k'}) \leq \epsilon$,

$$|(W_{j,k'}\sigma_{j,k'})^2 - (n_{j,k})^2| \leq c \left\{ (n_{j,k})^2 \left(1 - \mathcal{B}^{\frac{1-M}{1+M}j}\right)^{-1} + \varrho_j \left(\epsilon + \mathcal{B}^{\frac{1-M}{1+M}j}\right) \right\}.$$

The proof of Lemma 9 is postponed to the end of the Section. For

$$M \geq \frac{1+\kappa}{1-\kappa} \tag{18d}$$

and using Lemma 9 with $\epsilon = \frac{1}{2}\mathcal{B}^{-\kappa j}$ we get

$$
\begin{aligned}
A_4 &\lesssim \left(1 - \mathcal{B}^{\frac{1-M}{1+M}j}\right)^{-2} A_3 \\
&\quad + \varrho^2 \left(\frac{1}{2}\mathcal{B}^{-\kappa j} + \mathcal{B}^{\frac{1-M}{1+M}j}\right)^2 \sum_{\substack{k,k' \in \mathcal{K} \\ d(\xi_k, \xi_{k'}) > \frac{1}{2}\mathcal{B}^{-\kappa j}}} \frac{\sum_{p: d(\xi_k, \xi_p) \leq \frac{1}{2}\mathcal{B}^{-\kappa j}} \lambda_p |\psi_k(\xi_p)\psi_{k'}(\xi_p)|^2}{(C + n_k^2)^2(C + n_{k'}^2)^2} \\
&\lesssim A_3 + \varrho^2 \left(\frac{1}{2}\mathcal{B}^{-\kappa j} + \mathcal{B}^{\frac{1-M}{1+M}j}\right)^2 \sum_{k,k' \in \mathcal{K}} \frac{1}{(C+n_k^2)^2(C+n_{k'}^2)^2} \mathcal{B}^{4j}\mathcal{B}^{-j(1-\kappa)2M}\mathcal{B}^{-2\kappa j} \\
&\lesssim A_3 \left(1 + \varrho^2 \mathcal{B}^{-2\kappa j} \frac{1}{C^2}\right)
\end{aligned}
$$

as $\frac{1-M}{1+M} < -\kappa$. Note now that under ii) and using (1c) and (1e), we have

$$\varrho_j \underset{j \to \infty}{=} O\left(\mathcal{B}^{(\kappa_2 + 3)j}\right) ;$$

and

$$\frac{A_4}{C^2 S^2} \lesssim \frac{A_3}{C^2 S^2} + \mathcal{B}^{2j(2 - \kappa_1 + \kappa_2 + 3 + \alpha - (1-\kappa)M)} \underset{j \to \infty}{\longrightarrow} 0$$

provided that

$$M > \frac{5 + \alpha - \kappa_1 + \kappa_2}{1 - \kappa} . \tag{18e}$$

It remains to see that their exists a $\kappa \in (0,1)$ satisfying simultaneously the conditions (18a) to (18e). This is the case if $M > \max\left\{\frac{\kappa_2 + \alpha}{\kappa_1 - 1} ; \frac{4 - \kappa_1}{\kappa_1 - 1} ; \frac{5 + \alpha - \kappa_1 + \kappa_2}{\kappa_1 - 1}\right\}$. It suffices to take $\kappa \in (2 - \kappa_1, 1)$ sufficiently near to $2 - \kappa_1$ (recall that $\kappa_1 > 1$). Finally, remark that any $\kappa_1$ satisfying i is bounded by 2, then the above condition on $M$ reduces to the one mentioned in the statement of the Theorem. This concludes the proof of (14) and Theorem 2.

**Proof of Lemma 9**

First, for any $\delta > 0$, using (1b), (1e), (8), $0 \leq W \leq 1$ and the definition of $\varrho$

$$\begin{aligned}
|W_k^2 \sigma_k^2 - n_k^2| &= \left| W_k^2 \sigma_k^2 - \frac{1}{\lambda_k} \sum_{p=1}^N \lambda_p^2 W_p^2 \sigma_p^2 \psi_k^2(\xi_p) \right| \\
&\leq \frac{1}{\lambda_k} \sum_{p=1}^N \lambda_p^2 \left| W_k^2 \sigma_k^2 - W_p^2 \sigma_p^2 \right| \psi_k^2(\xi_p) + W_k^2 \sigma_k^2 \left| 1 - \frac{1}{\lambda_k} \sum_{p=1}^N W_p^2 \lambda_p^2 \psi_k^2(\xi_p) \right| \\
&\lesssim \sum_{p:d(\xi_p,\xi_k)\leq\delta} \lambda_p |W_k^2 \sigma_k^2 - W_p^2 \sigma_p^2| \psi_k^2(\xi_p) + \sum_{p:d(\xi_p,\xi_k)>\delta} \lambda_p |W_k^2 \sigma_k^2 - W_p^2 \sigma_p^2| \psi_k^2(\xi_p) \\
&\quad + W_k^2 \sigma_k^2 \left( 1 + \sum_{p=1}^N \lambda_p \psi_k^2(\xi_p) \right) \\
&\lesssim \varrho\delta + \frac{\mathcal{B}^j}{(1+\mathcal{B}^j\delta)^M} \sum_{p:d(\xi_p,\xi_k)>\delta} \lambda_p \left| W_k^2 \sigma_k^2 - W_p^2 \sigma_p^2 \right| + W_k^2 \sigma_k^2 \\
&\lesssim \varrho \left( \delta + \frac{\mathcal{B}^j}{(1+\mathcal{B}^j\delta)^M} \right) + W_k^2 \sigma_k^2 \quad . \tag{19}
\end{aligned}$$

Second

$$\begin{aligned}
n_k^2 &\asymp \sum_{p=1}^N \lambda_p W_p^2 \sigma_p^2 \psi_k^2(\xi_p) = W_k^2 \sigma_k^2 \sum_{p=1}^N \lambda_p \psi_k^2(\xi_p) - \sum_{p=1}^N \lambda_p \left( W_k^2 \sigma_k^2 - W_p^2 \sigma_p^2 \right) \psi_k^2(\xi_p) \\
&= W_k^2 \sigma_k^2 \sum_{p=1}^N \lambda_p \psi_k^2(\xi_p) - \sum_{p:d(\xi_p,\xi_k)\leq\delta} \left( W_k^2 \sigma_k^2 - W_p^2 \sigma_p^2 \right) \psi_k^2(\xi_p) - \sum_{p:d(\xi_p,\xi_k)>\delta} \left( W_k^2 \sigma_k^2 - W_p^2 \sigma_p^2 \right) \psi_k^2(\xi_p) \\
&\gtrsim W_k^2 \sigma_k^2 \left( 1 - \frac{\mathcal{B}^j}{(1+\mathcal{B}^j\delta)^M} \right) - \varrho \left( \delta + \frac{\mathcal{B}^j}{(1+\mathcal{B}^j\delta)^M} \right) \quad . \tag{20}
\end{aligned}$$

Combining Eqs. (19) and (20), one gets, for any $k'$ with $d(\xi_k, \xi_{k'}) \leq \epsilon$

$$\begin{aligned}
|W_{k'}^2 \sigma_{k'}^2 - n_k^2| &\leq |W_{k'}^2 \sigma_{k'}^2 - W_k^2 \sigma_k^2| + |W_k^2 \sigma_k^2 - n_k^2| \\
&\lesssim \varrho\epsilon + \varrho \left( \delta + \frac{\mathcal{B}^j}{(1+\mathcal{B}^j\delta)^M} \right) \\
&\quad + \left( n_k^2 + \varrho \left( \delta + \frac{\mathcal{B}^j}{(1+\mathcal{B}^j\delta)^M} \right) \right) \left( 1 - \frac{\mathcal{B}^j}{(1+\mathcal{B}^j\delta)^M} \right)^{-1} \\
&= n_k^2 \left( 1 - \frac{\mathcal{B}^j}{(1+\mathcal{B}^j\delta)^M} \right)^{-1} \\
&\quad + \varrho \left[ \epsilon + \left( \delta + \frac{\mathcal{B}^j}{(1+\mathcal{B}^j\delta)^M} \right) \left( 1 + \left( 1 - \frac{\mathcal{B}^j}{(1+\mathcal{B}^j\delta)^M} \right)^{-1} \right) \right] \quad .
\end{aligned}$$

Particularizing this last expression with $\delta = \mathcal{B}^{\frac{1-M}{1+M}j}$ leads to Lemma 9, since we have, for $M \geq 3$, $j \geq 1$ and this $\delta$: $\frac{\mathcal{B}^j}{(1+\mathcal{B}^j\delta)^M} < \frac{1}{2}$ .

**References**

[1] P. Baldi, G. Kerkyacharian, D. Marinucci, and D. Picard. Asymptotics for spherical needlets. *Ann. Statist.*, 2008. (In press).